%% file: main.tex
\documentclass[review]{elsarticle}
\usepackage{amsmath,amssymb,amsfonts}
\usepackage{algorithmic}
\usepackage{graphicx}
\usepackage{textcomp}
\usepackage{xcolor}

\usepackage{multirow}
\usepackage{url}
\usepackage[utf8]{inputenc}
\usepackage{verbatim}
\usepackage{nomencl}
\usepackage{tabularx,booktabs}
\usepackage[linesnumbered,ruled,vlined]{algorithm2e}

\def\BibTeX{{\rm B\kern-.05em{\sc i\kern-.025em b}\kern-.08em
    T\kern-.1667em\lower.7ex\hbox{E}\kern-.125emX}}

\usepackage{lineno,hyperref}
\modulolinenumbers[1]

\usepackage{scalerel}
\usepackage{tikz}
\usetikzlibrary{svg.path}

\definecolor{orcidlogocol}{HTML}{A6CE39}
\tikzset{
  orcidlogo/.pic={
    \fill[orcidlogocol] svg{M256,128c0,70.7-57.3,128-128,128C57.3,256,0,198.7,0,128C0,57.3,57.3,0,128,0C198.7,0,256,57.3,256,128z};
    \fill[white] svg{M86.3,186.2H70.9V79.1h15.4v48.4V186.2z}
                 svg{M108.9,79.1h41.6c39.6,0,57,28.3,57,53.6c0,27.5-21.5,53.6-56.8,53.6h-41.8V79.1z M124.3,172.4h24.5c34.9,0,42.9-26.5,42.9-39.7c0-21.5-13.7-39.7-43.7-39.7h-23.7V172.4z}
                 svg{M88.7,56.8c0,5.5-4.5,10.1-10.1,10.1c-5.6,0-10.1-4.6-10.1-10.1c0-5.6,4.5-10.1,10.1-10.1C84.2,46.7,88.7,51.3,88.7,56.8z};
  }
}

\newcommand\orcidicon[1]{\href{https://orcid.org/#1}{\mbox{\scalerel*{
\begin{tikzpicture}[yscale=-1,transform shape]
\pic{orcidlogo};
\end{tikzpicture}
}{|}}}}

\usepackage{hyperref} 
\usepackage{cleveref}

\begin{document}

\renewcommand*{\today}{February 25, 2021}
\begin{frontmatter}

\title{Long-term Value of Flexibility from Flexible Assets in Building Operation}

\author{Kasper Emil Thorvaldsen\fnref{myfootnote}}\corref{mycorrespondingauthor}
\cortext[mycorrespondingauthor]{Corresponding author}
\ead{kasper.e.thorvaldsen@ntnu.no}
\author{Magnus Korpås\fnref{myfootnote}}
\author{Hossein Farahmand\fnref{myfootnote}}

\address{Department of Electric Power Engineering, Norwegian University of Science and Technology, 7491 Trondheim, Norway}

\begin{abstract}
In this work, we investigate how flexible assets within a residential building influence the long-term impact of operation. We use a measured-peak grid tariff (MPGT) that puts a cost on the highest single-hour peak import over the month. We apply a mathematical model of a Home Energy Management System (HEMS) together with Stochastic Dynamic Programming (SDP), which calculates the long-term impact of operating as a non-linear expected future cost curve (EFCC) from the end of the scheduling period to the start. The proposed model is applied to a case study for a Norwegian building with smart control of a battery energy storage system (BESS), Electric vehicle (EV) charging and space heating (SH). Each of the flexible assets are investigated individually with MPGT and for an energy-based grid tariff. The results showed that EV charging has the highest peak-power impact in the system, decreasing the total electricity cost by 14.6\% with MPGT when controllable compared to a reference case with passive charging. It is further shown how the EFCC helps achieve optimal timing and level of the peak demand, where it co-optimizes both real-time pricing and the MPGT.

\end{abstract}


\end{frontmatter}

\newpage

\input{chapters/Nomenclature}
\input{chapters/Introduction}

\input{chapters/Methodology_new}

\input{chapters/Case_study}

\input{chapters/Results_and_discussions}

\input{chapters/Conclusion}

\bibliographystyle{elsarticle-num}
\bibliography{ref_no_url}

\newpage
\pagebreak

\end{document}

%% file: chapters/Nomenclature.tex
\makenomenclature

\renewcommand{\nomgroup}[1]{%
\ifthenelse{\equal{#1}{A}}{\item[\textbf{SDP sets}]}{%
\ifthenelse{\equal{#1}{B}}{\item[\textbf{Index sets}]}{%
\ifthenelse{\equal{#1}{C}}{\item[\textbf{Parameters}]}{%
\ifthenelse{\equal{#1}{D}}{\item[\textbf{Decision variables}]}{%
\ifthenelse{\equal{#1}{E}}{\item[\textbf{Stochastic variables}]}{%
\ifthenelse{\equal{#1}{O}}{\item[\textbf{Other}]}{}}}}}}
}

\nomenclature[B]{$\mathcal{S}_{g}$}{set of state variables}
\nomenclature[BA]{$\mathcal{G}$}{set of days within the month}

\nomenclature[B]{$\mathcal{T}$}{Set of time steps within the day}

\nomenclature[E]{$\mathcal{C}^{spot}_t$}{Electricity spot price in time step t [$\frac{EUR}{kWh}$]}
\nomenclature[E]{$D_{t}^{el}$}{Consumer-specific load in time step t [$kWh$]}
\nomenclature[E]{$\delta_t^{EV}$}{EV connected to building $\{0,1\}$}
\nomenclature[E]{$I^{Irr}_t$}{Solar irradiation at building in time step t[$\frac{kWh}{m^2}$]}
\nomenclature[E]{$T^{out}_t$}{Outdoor temperature in time step t [$^\circ C$]}


\nomenclature[D]{$y_t^{imp}, y_t^{exp}$}{Energy imported/exported to household [$\frac{kWh}{h}$]}
\nomenclature[D]{$p^{imp}$}{Peak of imported energy [$\frac{kWh}{h}$]}
\nomenclature[D]{$\gamma$}{SOS-2 variables for the Expected future cost curve}
\nomenclature[D]{$\alpha_{p^{imp},s^e_{t+1}}^{future}$}{Expected future cost from peak power [$\frac{EUR}{\frac{kWh}{h}}$]}

\nomenclature[D]{$E_t^{EV}$}{State of charge for EV for time step t [$kWh$]}
\nomenclature[D]{$y_t^{EV,ch}$}{Input power to EV for time step t [$\frac{kWh}{h}$]}

\nomenclature[D]{$E_t^{B}$}{State of charge for Battery for time step t [$kWh$]}
\nomenclature[D]{$y_t^{B,ch},y_t^{B,dch}$}{Power to/from the battery for time step t [$\frac{kWh}{h}$]}

\nomenclature[D]{$y_t^{PV}$}{Power produced from PV system [$\frac{kWh}{h}$]}

\nomenclature[D]{$T_t^{in},T_t^{e}$}{Interior and building envelope temperature [$^\circ C$]}


\nomenclature[D]{$q_t^{sh}$}{Power usage for space heating [$\frac{kWh}{h}$]}


\nomenclature[C]{$\mathcal{C}_n^{imp}$}{Expected future cost for point $n$ [$EUR$]}
\nomenclature[C]{$\mathcal{C}^{grid}$}{DSO energy tariff for imported energy [$\frac{EUR}{kWh}$]}
\nomenclature[C]{$\mathcal{C}^{peak}$}{DSO capacity-based tariff for highest peak import [$\frac{EUR}{\frac{kWh}{h}}$]}
\nomenclature[C]{$P^{imp,max}$}{Maximum power import to building [$\frac{kWh}{h}$]}
\nomenclature[C]{$P_n^{imp}$}{Peak power at point $n$ [$\frac{kWh}{h}$]}
\nomenclature[C]{$P_0^{imp}$}{Initial peak power [$\frac{kWh}{h}$]}
\nomenclature[C]{$VAT$}{Value added tax for purchase of electricity [$p.u$]}

\nomenclature[C]{$E^{EV,Cap}$}{EV storage capacity [$kWh$]}
\nomenclature[C]{$\eta^{EV}_{ch}$}{EV charging efficiency [$\%$]}
\nomenclature[C]{$\dot{E}^{Max}$}{Maximum EV charging capacity [$\frac{kWh}{h}$]}
\nomenclature[C]{$E^{EV,min},E^{EV,max}$}{Min/Max EV SoC capacity [$kWh$]}
\nomenclature[C]{$D^{EV}$}{EV discharge when not connected [$kWh$]}

\nomenclature[C]{$E^{B,Cap}$}{Battery storage capacity [$kWh$]}
\nomenclature[C]{$\eta^{B}_{dch},\eta^{B}_{ch}$}{Discharge/charge efficiency for battery [$\%$]}
\nomenclature[C]{$\dot{E}^{B,dch},\dot{E}^{B,ch}$}{Discharge/charge capacity for battery [$\frac{kWh}{h}$]}

\nomenclature[C]{$E^{B,min},E^{B,max}$}{Battery SoC limits [$kWh$]}

\nomenclature[C]{$A^{PV}$}{PV system area [$m^2$]}
\nomenclature[C]{$\eta^{PV}$}{Total efficiency for PV system [$\%$]}

\nomenclature[C]{$\dot{Q}^{sh}$}{Capacity for space heating radiator [$\frac{kWh}{h}$]}

\nomenclature[C]{$T_t^{in,min}, T_t^{in,max}$}{Lower/upper interior boundary [$^\circ C$]}

\nomenclature[C]{$R_{ie}, R_{eo}$}{The thermal resistance between the interior-building envelope and building envelope-outdoor area [$\frac{^\circ C}{kWh}$]}
\nomenclature[C]{$C_i, C_e$}{Heat capacity for interior and building envelope [$\frac{kWh}{^\circ C}$]}


\nomenclature[C]{$N_S$}{Number of nodes for stochastic variables}
\nomenclature[C]{$N_P$}{Number of discrete peak power values}

\printnomenclature

%% file: chapters/Introduction.tex
\section{Introduction} \label{Introduction}

With the roll-out of smart meters, it is possible to implement more dynamic price structures so the end-users can react to price signals on their own accord. This introduces the potential of participating in price-based Demand Response (DR) programs, which ideally should be able to represent the grid operators' actual system cost of operating the grid. However, it is important that the system cost is accurately represented in their programs. 

The electricity bill from the grid operator for the end-user is currently in a period of transition, where the cost is going from a passive volumetric charge cost to a combination of volumetric and capacity-based costs, in line with recommendations \cite{EDSOforSmartGrids2015AdaptingFuture}. The measured-peak grid tariff (MPGT) and subscribed capacity have been presented by the Norwegian Water Resources and Energy Directorate (NVE) as possible new grid tariffs in Norway \cite{Norgesvassdrags-ogenergidirektorat2017ForslagNettvirksomhet}. MPGT introduces a capacity-based tariff determined by the highest single-hour energy consumption over an hour for a given period. 
The MPGT is already implemented for end-users with a total yearly consumption rate over 100 MWh, for a monthly period \cite{AgderAE.no}. If such a grid tariff is implemented for smaller consumers, a Home Energy Management System (HEMS) can help consumers achieve a more cost-optimal utilization of local energy resources in response to such tariffs if they are capable of operating while considering the whole period. 


In this study, we focus on the impact of flexible assets in a HEMS, such as battery energy storage system (BESS), thermal flexibility from space heating (SH) and electric vehicles (EVs). In \cite{Zheng2015SmartShaving}, different BESS alternatives are investigated in an economic overview for an average residential consumer in the US with the overall goal of shaving electrical peaks under a Time-of-use demand tariff. A thermal storage tank was used in \cite{Zhang2013EfficientMicrogrid} for a building consisting of multiple smart homes with wind production and BESS to reduce peak electricity consumption from the grid. The work in \cite{Berge2016PerceivedBuildings} analyzed thermal comfort and temperature zoning in residential buildings with user feedback to analyze performance. Trying to charge an EV optimally given the uncertainty in driving pattern was investigated in \cite{Iversen2014OptimalProcess}, where the stochastic nature gave a noticeable impact on charging strategies. In \cite{Angenendt2019OptimizationCoupling}, the EV was used together with a HEMS and a charge-discharge management framework to charge the EV optimally while lowering the PV curtailment and reduce total residential operation cost. 

A HEMS will operate and control the energy input and output from the different flexible assets, adjusting the flow of energy based on what is deemed the optimal decision for the HEMS to consider. In most cases, this control of flexible assets is used to optimize the total cost of electricity within the period that the HEMS considers. However, the information that the HEMS consider for operation of the energy system, is subject to uncertainty. A literature review was conducted in \cite{Tian2018AAssessment} regarding uncertainty within building energy systems, showcasing how weather effects, the modeling of the building envelope, and occupant behavior are different kinds of input data uncertainty that should be accounted for. 
In \cite{Pallonetto2020OnBuildings}, an overview over DR-programs and their development was done, commenting on the importance of smart control systems for residential buildings to participate efficiently in these programs. In addition, different modeling techniques for a residential building and control algorithms for a HEMS were presented. However, they also raised the issue that appropriate tariff structures depend on the end-users capability of controlling their own consumption. 

Within the operation of a HEMS, multiple optimization methods and approaches have been presented in the research literature. A mixed-integer linear problem (MILP) definition was defined within a HEMS in \cite{Erdinc2014EconomicHouseholds}. The HEMS included photovoltaics (PV), BESS and EV with bi-directional power flow, and the model optimizes the system operation with dynamic pricing and peak power limits. Another approach was used for a smart building consisting of PV, heat pump, thermal storage and BESS in \cite{Kuboth2019EconomicSystems}, where the authors used a model predictive control (MPC) approach on a stochastic problem to optimize each of the flexible assets. A rolling horizon strategy was deployed in \cite{Palma-Behnke2013AStrategy} to operate an energy management system for a microgrid. The microgrid consisted of PV, wind turbines, a diesel generator and a BESS, together with demand-side management (DSM).

The authors in \cite{Salpakari2016OptimalPV} apply a deterministic dynamic programming (DP) model to analyze the cost-optimal control for a building with varying degree of PV installed. The application resulted in cost savings compared to a rule-based approach which maximizes PV self-consumption.
The work presented in \cite{bjarghov2017utilizing} developed a deterministic DP model that optimized the state-of-charge for either a BESS or EV battery over a year for a household with PV. For a power-based grid tariff, cost savings at 13.3\% with optimal battery control, or 16.6\% for optimal EV battery control were achieved.
The use of stochastic DP (SDP) in a HEMS has been investigated in \cite{Donadee2012StochasticVehicles}, where the authors utilized SDP to optimize both EV charging and frequency regulation bids given the expected future costs calculated by SDP. A smart building model is presented in \cite{Wu2016StochasticArray} that used SDP to optimize energy management for EV and PV with uncertainty in generation, consumption and EV availability. The approach found potential cost savings of almost 500\% for a Tesla Model S compared to no optimal control, and load shifting potential for interaction with the grid.
The work was further carried on in \cite{Wu2018StochasticStorage} where the authors analyzed several operating modes of the EV, such as V2G, V2H and G2V. The study found that by utilizing V2G in a grid with real-time pricing (RTP) and limited bi-directional grid capacity during certain times of the day, 75.5\% cost savings could be achieved compared to a case without a plug-in EV. Ref
\cite{Kim2011SchedulingUncertainty} investigated an energy consumption scheduling problem with uncertain future prices. The SDP algorithm was used to describe an optimal scheduling algorithm for non- and interruptible loads based on price thresholds.

Most of the work here using SDP considers a short-time horizon of up to a day, except for \cite{Salpakari2016OptimalPV,bjarghov2017utilizing}, both considering a year, but not indicating that their models could be used as a short-term operational model.
One of the most valuable takeaways from the use of SDP is that the scheduling period is decomposed into smaller decoupled segments, only coupled through expected future cost curves which depict the future consequence of altering the state variables.
Therefore, SDP can contribute on foreshadowing the expected impact of the system beyond the current scheduling horizon.  

Based on a review of the literature, the relationship between long-term and short-term value of energy flexibility has not yet been specified or investigated in full detail. To the best of the authors' knowledge, the inclusion of long-term price signals in building operational models has not previously been performed with the exception of \cite{EmilThorvaldsen2020RepresentingProgramming}, which includes our initial approach to proposing this idea. There, we applied SDP in a MPGT setting over the course of a month.

\begin{figure}[h!]
\centering
\includegraphics[width=1\textwidth]{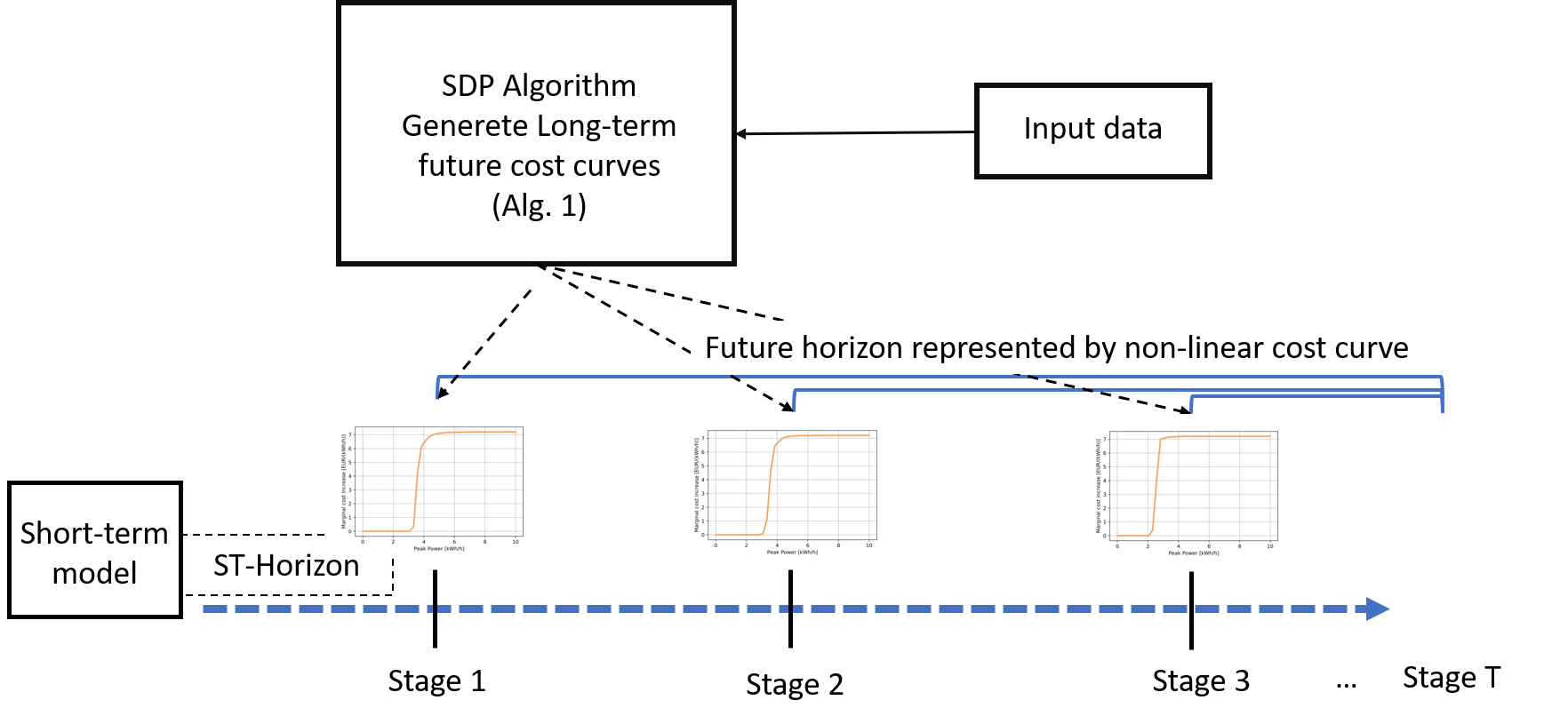}
\caption{Overview over the coupling of long-term future cost curves generated by the SDP framework in this work and in \cite{EmilThorvaldsen2020RepresentingProgramming}, and a short-term operational model. Alg. 1 will be further explained in Section \ref{Model}.}
\label{fig:SDP_framework}
\centering
\end{figure}

The methodology presented in \cite{EmilThorvaldsen2020RepresentingProgramming} decomposed one month (long-term horizon) into smaller daily stages (short-term horizon) and generated non-linear expected future cost curves (EFCC). The EFCC describes the future consequence of flexibility utilization within a HEMS at each stage, and can be given as an input to a detailed short-term operational model. The coupling of the proposed long-term operation to a short-term operating model is showcased in Fig. \ref{fig:SDP_framework}.


The SDP algorithm is used to generate EFCCs, showcasing the change in the expected cost beyond the decisions taken up to each stage. The SDP algorithm determines how a HEMS would react to a number of possible future scenarios, and calculates the cost-optimal decision backward, which then is put together into a weighted future cost curve and used in a short-term operational model. 
This approach is similar to hydropower scheduling, where long- and medium-term models are used to generate cost curves as showcased in Ref \cite{Gjelsvik1999AnUncertainty}, which are given as input to a short-term operational model \cite{Fodstad2015HydropowerMarkets}. This paper focuses on the generation of these future cost curves.



Expanding on the work presented in \cite{EmilThorvaldsen2020RepresentingProgramming}, the main contributions here are as follows:

\begin{itemize}
    \item We present a general SDP framework for optimal energy management of HEMS exposed to a monthly MPGT. The HEMS considers that SH, EV charging and a BESS can be utilized to keep the peak import at the cost-optimal level, influenced by the implications given by the EFCC calculated
    \item The output of the SDP framework is further analyzed for each individual flexible asset, to capture the different characteristics each asset contains in a long-term operational setting
    \item In a numerical case study based on real Norwegian conditions, the model is applied to two different electricity grid tariffs, and tested where each flexible asset can be controlled individually. The analysis showcases their impact and value of flexibility based on the peak import level, and each asset and the different schemes are compared against each other 
\end{itemize}

The paper progresses as follows. Section \ref{Model} introduces the SDP framework methodology of the HEMS. The case study used for the analysis is presented in Section \ref{Case_study}, while the results and discussion are found in Section \ref{discussion}. The conclusion and future work are given in Section \ref{Conclusion}.

%% file: chapters/Methodology_new.tex
\section{Methodology} \label{Model}

The overall objective of a HEMS is to minimize the expected total electricity cost from interaction with the distributional grid, consisting of electricity purchase and grid tariff payments. The scheduling horizon depends on the value of utilizing load shifting over longer periods, and the duration of any long-term price signals that are included. This work considers the long-term price signal from the MPGT, paid at the end of the horizon. Thus, the objective will be to find the optimal operation to obtain the minimum expected operating cost for the HEMS, as shown in Eq. (\ref{eq:sub_period}):

\begin{align}
    &  min \mathbb{E} \{ \sum_{t= 1}^{T_{sub}} [\mathcal{C}_t^{spot}\cdot (y_t^{imp} - y_t^{exp}) + \mathcal{C}^{grid}\cdot y_t^{imp}] + \Phi(p^{imp})  \}  \label{eq:sub_period}
\end{align}


We deal with this issue in a similar fashion to hydropower scheduling \cite{Helseth2013APower}, by assuming that we only want to optimize the HEMS for a pre-determined period between $t=1$ and $T_{sub}$ of the total horizon. For this pre-determined period to be capable of operating optimally and still consider the cost paid at the end of the horizon, the future impact must be included. Here, the function $\Phi (p^{imp})$ represents the cost for the peak import beyond the present period, into the future. 
Thus, the peak import $p^{imp}$ is coupled in time, giving a dynamic connection of the optimization problem. Hence, the formulation represents a multi-stage multi-scenario optimization problem, and we can apply decomposition techniques to simplify the problem. We  utilize an SDP approach for multi-stage in a backwards procedure \cite{Birge2011IntroductionOptimization}. The function $\Phi (p^{imp})$ from Eq. (\ref{eq:sub_period}) is represented as a piecewise-linear future cost curve calculated through an SDP algorithm.

The original problem is decomposed into several smaller single-stage deterministic problems to decrease the computational complexity.
The decomposition is performed through a set of state variables $\mathcal{S}_g$, which consists of all information that is carried over between the decision stages $g-1$ to $g$. This set consists of two unique subsets; subset $\mathcal{S}_{S,g}$, which consists of realized values for the stochastic variables for each decision stage $g$, and subset $\mathcal{S}_{P,g}$, which contains state variables in the optimization problem. 
Based on the set $\mathcal{S}_g$, a decomposed decision problem for the HEMS to solve is defined based on the state ${s_g^s,s_g^p} \in \mathcal{S}_g$ for a decision stage $g$, which is given by the current scheduling for that decision stage, and the weighted impact of the future cost for all scenarios. The decomposed decision problem is presented in Section \ref{Optimzation_model}.

To perform valid coupling between the stages, the conditions at the end of one stage to the next corresponding stage must be identical. If not ensured, the transition is not feasible and will result in inaccurate results. In a building, this would be connected to the energy levels and connectivity of certain shiftable units. If these variables are not tied to the state variable subset $\mathcal{S}_{P,g}$, such that the future impact of any changes are included, then it is possible to lock their start/end values between stages to make the transition feasible. This work has a set start/end condition for each stage for the flexible assets, described in Section \ref{sol_strat}.

The stages are considered to be decoupled from each other, while the scenarios between stages are tied together through a transition probability, represented as a Markov decision process (MDP). To enable a method for generating future cost curves in a backwards DP strategy, the problem is represented as a Markov decision  \cite{Process1957DOC}. 
The scenarios generated specify the uncertain parameters in the system, referred to as stochastic variables.
For a specific stage and scenario, the stochastic variables are realized as input for the HEMS. The future impact beyond the stage is affected by transitioning the scenarios forward and their probabilities. 
The transition probabilities are assumed to have a Markov property, which specifies the stochastic process is memoryless \cite{Gudivada2015BigApplications}. 
The transitioning scenarios and their impact in the future are represented as the future cost $\Phi(p^{imp})$, making it an MDP. Therefore, the impact of the future scenarios in the MDP can be represented in the present scenario as an EFCC. By combining the different discrete scenarios $\mathcal{N}_{S}$ in $\mathcal{S}_{S,g+1}$, the EFCC illustrates the weighted cost-based impact towards the future.






\input{chapters/Building_model}

\subsection{Solution Strategy} \label{sol_strat}

\begin{algorithm}
\For{$g = \mathcal{G}, \mathcal{G}-1,..,1$}
 {
    \For {$n \in \mathcal{N}_{P}$}
    {
    $P_0^{imp} \leftarrow P_n^{imp}$\\
        \For{$s_g^s \in \mathcal{N}_{S}$}
        {
            $\{ \mathcal{C}_t^{spot}, D_t^{El}, \delta_t^{EV}, I_t^{Irr},T_t^{out}\} \leftarrow \Gamma(g,{s_g^s}) $
            
            $\mathcal{C}^{imp}_{i} \leftarrow \Phi(i,s_{g}^s,g+1)$ for $i = 1..\mathcal{N}_{P}$
            
            $\mathcal{C}_{s_g^s,n} \leftarrow Optimize (\ref{eq:Obj_function}) - (\ref{eq:heat})$ 
            
        }
        
        \For{$s_{g-1}^s \in \mathcal{N}_{S}$}
        {
        $\Phi(n,s_{g-1}^s,g) = \sum_{s^s_{g} = 1}^{\mathcal{N}_{S}} \mathcal{C}_{s_g^s,n}\cdot \rho(g,s^s_{g}|s_{g-1}^s)$ 
        
        }
        
    }
    
 }

\caption{The SDP algorithm}
\label{Alg_1}
\end{algorithm}

The complete optimization problem is solved through the SDP algorithm solution strategy as shown in Algorithm \ref{Alg_1}. This process will generate EFCCs for every stage of the overall problem.

The SDP algorithm initiates at the last stage of the horizon, and computes the overall cost of the decomposed decision problem in Section \ref{Optimzation_model} for the number of discrete points of the state variable $n \in \mathcal{N}_P$, and the number of scenarios $s_g^s \in \mathcal{N}_S$. 
For each scenario, we realize the stochastic variables with values from $\Gamma$ in line 5 as input into the specific problem, and in line 6 the EFCC for the future scheduling day $g+1$ is realized for each discrete point in $\mathcal{C}_i^{imp}$. $\Gamma$ contains the realized stochastic variables based on the scheduling day $g$ and scenario $s_g^s$. The EFCC used are the results from the previous stage $g+1$, and for the initial case $g = \mathcal{G}$, the end cost for peak import by the MPGT is used.

In line 7, the optimization problem for the HEMS is solved in order to find the objective function value, which is affected by both the cost for operating within the stage and the resulting EFC beyond the period. 
To enable a feasible transition between stages, $T_t^{in}$, $T_t^e$, $E_t^{EV}$ and $E_t^{B}$ have a start/end condition the optimization problem must hold. There is a high penalty cost included for the SH variables if this condition cannot be met, but this penalty cost is not included in the generation of EFCCs as it has no further influence on the decision-making.

The result from the problem is then used to derive the EFCC $\Phi(n,s_{g-1}^e,g)$ for $n \in \mathcal{N}_P$. The calculation of the EFCC is performed in line 8-9, where the EFC for a specific state variable point is derived. For a state variable value $n$, the EFC representation for stage $g$ is calculated as the weighted future cost value for all scenarios that can occur in this stage, to be used in the previous stage $g-1$. The EFC is calculated for each scenario that occurs in stage $g-1$ through the loop in line 8, and the weight of each EFC is based on the corresponding transition probabilities $\rho(g,s^s_{g}|s_{g-1}^s)$, which depends on the probability of arriving at scenario $s_g^s$ when originally at scenario $s_{g-1}^s$. This process is performed for all state variables until the whole EFCC has been calculated for all scenarios $s_{g-1}^s \in \mathcal{N}_S$. These results are then used as the basis for stage $g-1$, until we have arrived at the first stage and have derived EFCCs for all stages and scenarios.

%% file: chapters/Building_model.tex
\subsection{Decomposed decision problem} \label{Optimzation_model}

\begin{figure}[h!]
\centering
\includegraphics[width=0.6\textwidth]{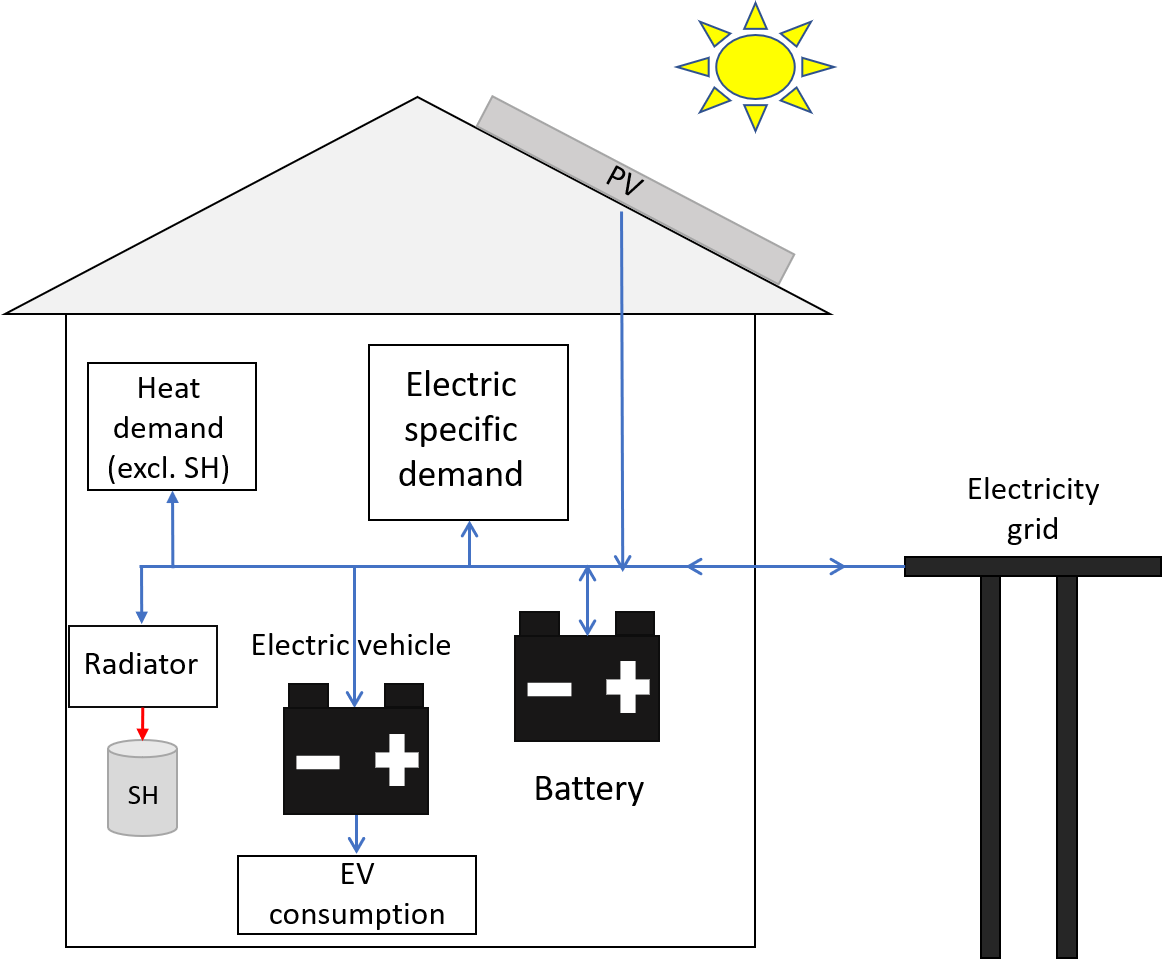}
\caption{Overview of building structure and energy system.}
\label{fig:Building}
\end{figure}

The decomposed decision problem is formulated as an optimization model for a given stage $g$, scenario $s_g^s$, and initial peak import power $P_{0}^{imp}$. The optimization model described in this subsection consists of a residential building connected to the power grid with bi-directional power flow capability, as illustrated in Fig. \ref{fig:Building}, operated through a HEMS. 
The electric-specific demand and heat demand from a water tank must be met at all times and these are treated as non-shiftable loads $D_t^{El}$. The smart control covers SH, BESS control, EV charging, and a roof-mounted photovoltaic (PV) system.

\subsubsection{Objective function}

The objective function in (\ref{eq:Obj_function}) minimizes the total electricity cost for the end-user. The electricity cost is the cost or benefit of importing/exporting electricity from the grid, respectively, and the future long-term cost for operating beyond this stage based on the highest single-hour peak import power $p^{imp}$, based on the EFCC. The objective function will co-optimize both short-term and long-term implications, where it balances the marginal cost increase of short-term operation versus the marginal cost savings in the long-term when lowering the peak import. 

\begin{align}
    &  min \{ \sum_{t\in \mathcal{T}} [\mathcal{C}_t^{spot}\cdot (y_t^{imp} - y_t^{exp}) + \mathcal{C}^{grid}\cdot y_t^{imp}] + \alpha_{p^{imp},s_{g+1}^{s}}^{future}  \}  \label{eq:Obj_function}
\end{align}  


\subsubsection{Energy balance}
The electric energy balance of the house is given in (\ref{eq:Energy_balance}). The exchange between the grid, together with production from the PV system and interaction from the BESS, must cover the needed load for the building.

\begin{align}
    & y_t^{imp} - y_t^{exp} + y_t^{PV}  + y_t^{B,dch} = \nonumber \\
    & D_t^{El} + y_t^{EV,ch} + q_t^{sh} + y_t^{B,ch} \quad \forall t \label{eq:Energy_balance}
\end{align}

\subsubsection{Expected Future Cost Curve} \label{EFCC_expl}

The EFC for this problem is depicted within (\ref{eq:Memory_power}) to (\ref{eq:Non_negativity_sos2}). The highest amount of power that is imported to the building is denoted by $p^{imp}$, which is bounded by the highest peak within the decision stage (\ref{eq:hourly_power}) and the initial value from earlier periods (\ref{eq:Memory_power}). The peak import power achieved at the end is used to set the EFC included in the objective function, which consists of discretized values of
$P_n^{imp} \; n \in \mathcal{N}_{P}$ represented through SOS-2 variables 
\cite{Beale1976GLOBALSETS}. 


\begin{subequations}
    \begin{flalign}
        & p^{imp} \geq P_0^{imp} \label{eq:Memory_power}\\
        & p^{imp} \geq y_t^{imp} \quad \forall t \label{eq:hourly_power}\\
        & \alpha_{p^{imp},s_{g+1}^{s}}^{future} =   \sum_{n \in \mathcal{N}_{P}} \gamma_n \cdot \mathcal{C}_{n}^{imp} \label{eq:SOS2_decision} \\
        & p^{imp} = \sum_{n \in \mathcal{N}_{P}} \gamma_n \cdot P_n^{imp} \label{eq:SOS2_binding} \\
        & \sum_{n \in \mathcal{N}_{P}} \gamma_n = 1 \label{eq:SOS2_var} \\
        & \gamma_n \geq 0 \quad \forall n , \: SOS {\text -} 2 \label{eq:Non_negativity_sos2}
    \end{flalign}
\end{subequations}

\subsubsection{Electric vehicle}

The behavior of the EV is formulated in (\ref{eq:EV_soc_balance}) to (\ref{eq:EV_soc_lim}). The EV is modeled as a uni-directional battery that can be charged at a continuous rate, with an availability pattern based on the stochastic variable $\delta^{EV}_t$. If the EV is not available at the building, the EV cannot be charged and a constant load discharges the battery to simulate consumption from driving. The EV must stay within a specific range in its state-of-charge (SoC) in (\ref{eq:EV_soc_lim}), whereas the boundary is time-dependent to include travelling preferences.

\begin{subequations}
    \begin{align}
    & E_t^{EV}-E_{t-1}^{EV} = y_t^{EV,ch}\eta^{EV}_{ch}\delta_t^{EV} \nonumber \\
    & -D^{EV} (1-\delta_t^{EV}) \quad  \quad \forall t  \label{eq:EV_soc_balance} \\
    & 0 \leq y_t^{EV,ch} \leq \dot{E}^{Max} \quad \forall t \label{eq:EV_charge_lim} \\
    & E_t^{EV,min} \leq E_t^{EV} \leq E_t^{EV,max} \quad \forall t \label{eq:EV_soc_lim} 
    \end{align}
\end{subequations}

\subsubsection{Battery energy storage system}

A bi-directional BESS is available within the building with the characteristics shown in (\ref{eq:BAT_balance}) to (\ref{eq:Bat_soc_lim}). The battery can be discharged and charged at a continuous rate, with limitations on power capacity and a storage capacity range to ensure optimal operation without risk of damage.

\begin{subequations}
    \begin{align}
    & E_t^{B} - E_{t-1}^{B} = y_t^{B,ch}\eta^{B}_{ch} - \frac{y_t^{B,dch}}{\eta^{b}_{dch}} \; \quad \forall t \label{eq:BAT_balance} \\
    & 0 \leq y_t^{B,ch}\eta^{B}_{ch} \leq \dot{E}^{B,ch} \quad \forall t \in \mathcal{T} \label{eq:Bat_charge_lim} \\
    & 0 \leq y_t^{B,dch} \leq \dot{E}^{B,dch} \quad \forall t \in \mathcal{T} \label{eq:Bat_discharge_lim} \\
    & E^{B,min} \leq E_t^{B} \leq E^{B,max} \quad \forall t \label{eq:Bat_soc_lim}
    \end{align}
\end{subequations}


\subsubsection{Photovoltaic system}

A roof-mounted PV system is connected to the electrical system through a controllable system. 
The HEMS is assumed to change the power output in a similar fashion to the work presented in \cite{Jain2017AnSystem}.

\begin{align}
    &  0 \leq y_t^{PV} \leq A^{PV} \cdot \eta^{PV}\cdot I_t^{Irr} \quad \forall t \in \mathcal{T} \label{eq:PV_prod} 
\end{align}

\subsubsection{Space heating}

All considerations regarding heating of the building are formulated in (\ref{eq:Heater}) to (\ref{eq:Envelope_balance}). The building has an electric radiator for SH that can be operated continuously. The heat dynamics in the building are shown as a grey-box model, in which the physical behavior is formulated using linear state-space models \cite{Sonderegger1978DynamicParameters, Bacher2011IdentifyingBuildings}. The dynamics between the interior temperature and the outdoor temperature can be captured alongside disturbances as heat input, which will include the impact of time-dependent temperature deviations. Thus, the heat system can be represented through an RC-network model. In an RC-network model, resistors (R) represents thermal resistance between measuring points, capacitors (C) capture the heat capacity of the measuring point, and ${q}^{sh}$ is the heat flux from heat sources. In addition, the outdoor temperature impact is included as a voltage source ($T^{out}$).

\begin{figure}[h!]
\centering
\includegraphics[width=0.8\textwidth]{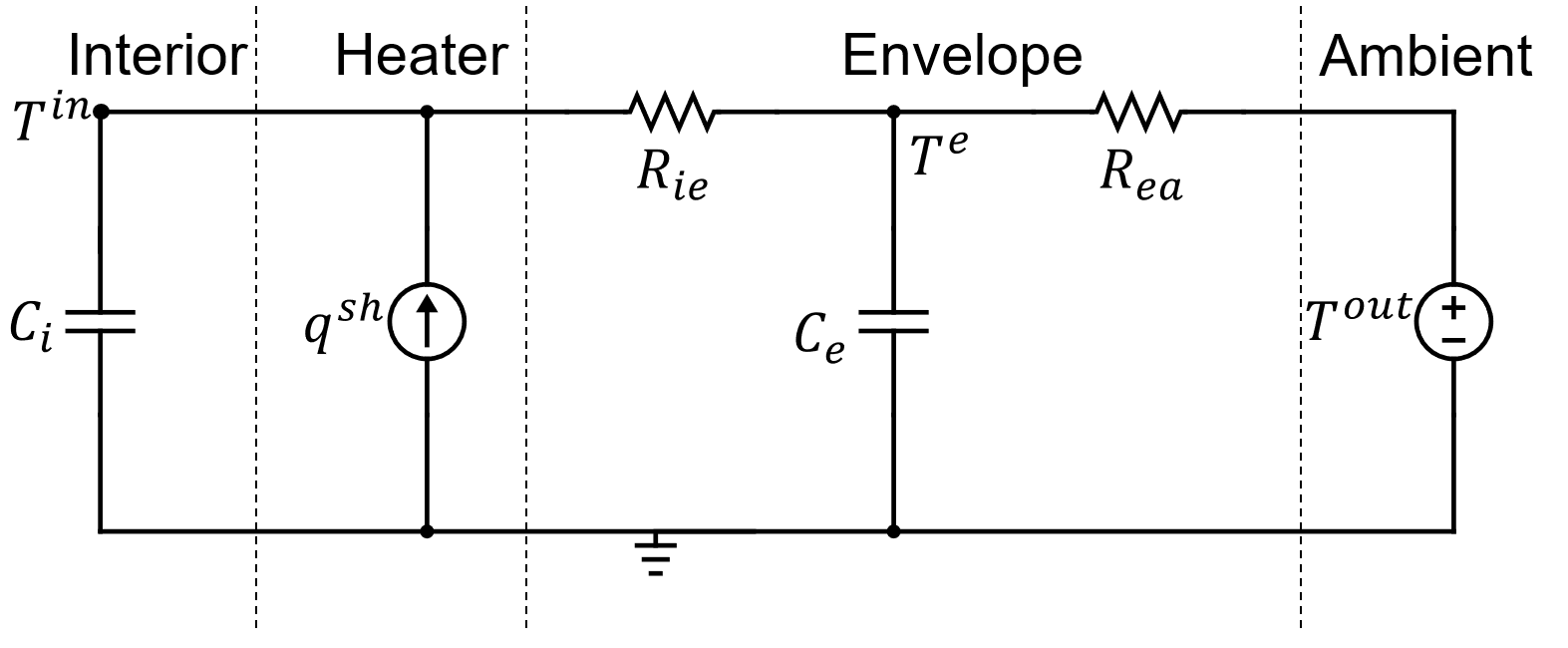}
\caption{RC-network for the SH dynamics given in equations (\ref{eq:interior_balance}) to (\ref{eq:Envelope_balance}).}
\label{fig:RC-network}
\end{figure}

The RC-network layout depends on the number of zones and inputs into each existing zone \cite{Bacher2011IdentifyingBuildings}. 
This optimization problem utilizes a 2R2C model, which divides the system into three zones: the interior or indoor of the building, the envelope acting as the physical separator between the interior and outdoors, and the outdoor area. The layout is shown in Fig. \ref{fig:RC-network}. The interior temperature is measured by the control system for the building, which can utilize the electric heater to regulate the temperature in response to impact from the envelope and outdoor temperature.

\begin{subequations}
    \begin{align}
    & 0 \leq q_t^{sh} \leq Q^{sh} \quad \forall t \label{eq:Heater} \\
    & T_t^{in,min} \leq T_t^{in} \leq T_t^{in,max} \quad \forall t \label{eq:Heat_lim_empty} \\
    & T^{in}_{t} - T^{in}_{t-1} = \frac {1}{R_{ie}C_i} [T^{e}_{t-1} - T^{in}_{t-1}] + \frac{1}{C_i}q_{t}^{sh} \quad \forall t  \label{eq:interior_balance} \\
    & T^{e}_{t} - T^{e}_{t-1} = \frac {1}{R_{ie}C_e}[T^{in}_{t-1} - T^{e}_{t-1}] \nonumber \\
    & + \frac{1}{R_{eo}C_i} (T_{t-1}^{out}-T_{t-1}^e) \quad \forall t \label{eq:Envelope_balance} 
    \end{align}
    \label{eq:heat}
\end{subequations}

%% file: chapters/Case_study.tex
\section{Case Study} \label{Case_study}

The presented model has been applied to a realistic case study of a Norwegian building in which the presented MPGT is included in the electricity bill. The building is a single-family house (SFH) placed in the south-eastern part of Norway, and the HEMS controls the different flexible assets available. The analysis is for January 2017 with hourly time resolution per day, and the stochastic variables consist of historical data or synthetic data from supporting literature. 


\subsection{Building structure}

The PV system on the roof consists of 4.65 kW installed capacity, connected through an MPP inverter with a combined constant conversion and MPP efficiency at 95\% \cite{Valentini2008PVEvaluation}.

The inelastic demand originates from two sources: The electric-specific electricity consumption from users, and non-flexible domestic hot water (DHW) consumption covered by a water tank. The data for the electric-specific electricity are obtained from the Distributional Operator (DSO) Ringerikskraft from January 2017 \cite{Ringerikskraft}. 
The DHW-consumption profile is based on measurement of 49 water heaters at Norwegian households through the "Electric Demand Knowledge - ElDek" \footnote{https://www.sintef.no/prosjekter/eldek-electricity-demand-knowledge/} research project by SINTEF Energy Research \cite{SINTEFEnergyResearchKMBproject190780/S60ElDeKKnowledge}. The RTP electricity prices from bidding zone NO1 in Nordpool for year 2017 are used \cite{SeeYou.b}. 

\subsubsection{Heat dynamics}
The layout in Fig. \ref{fig:RC-network} represents the heat dynamics of the building, and is based on observed values from the Living Lab building built by Zero Emission Building (FME ZEB)\footnote{www.fmezen.no/} and NTNU \cite{LivingLab,Vogler-Finck2019InverseHouse}. The Living Lab is a pilot project used to study various technologies and design strategies with the overall goal of reaching the zero emission target and analyzing thermo-physical properties \cite{FinocchiaroTheFacility}. 
The space heating is performed through a 3 kW radiator, which can operate continuously. The default temperature boundary that the HEMS uses is a range of 20-24 $^\circ C$, which from reference \cite{Berge2016PerceivedBuildings} is the threshold end-users find comfortable.

\subsubsection{Electric Vehicle}
A 24 kWh EV is selected for this study, with an operational range between 20-90\% of total capacity at all times, with a range between 60-90\% when departing to prevent range anxiety. The EV consumes electricity from the battery when it has departed to simulate driving. 


The EV consumption rate over a weekday has been simplified as a deterministic input. Based on \cite{LakshmananValueScenario}, the mean driving distance of 52 km has been used, under the assumption that the EV consumes 18 kWh/100km, which puts the hourly discharge rate at 1.02 kWh/h for $D^{EV}$ with a 9-hour offline timeframe. The EV is assumed to always be connected during stage transition, and as the optimization model in Section \ref{Optimzation_model} has deterministic info for each stage and scenario analyzed, the HEMS can charge the EV so it ensures enough SoC for the trip.




\subsubsection{Battery energy storage system}
The BESS installed in this system is based on a battery from SonnenBatterie \cite{SonnenBatterie} with a rated power input/output of 2.5 kW measured at the output of the inverter. The tolerated SoC is set at between 10-100\% SoC, and a round-trip efficiency of 85\% in line with efficiency settings from \cite{Jafari2020PowerVariability}. Any cost or performance associated with battery degradation is left out of this analysis.

\subsubsection{Initial conditions}
As mentioned in Section \ref{sol_strat}, the following variables have been given a start/end value to enable a feasible stage transition: $T_0^{in} = 22$ $^\circ C$, $T_0^e = 20$ $^\circ C$, $E_0^{EV} = 14.4$ kWh, $E_0^{B} = 2.5 $ kWh. 

\subsubsection{Grid tariff structure}

The grid tariff structure consists of multiple layers of payment. A conversion rate of 1 EUR = 10 NOK has been used for this work. The first is a fixed consumer cost given as a volumetric cost at $0.024 \frac{EUR}{kWh}$ in 2017. The rest contains the cost provided by the DSO, which depends on the tariff scheme. In 2017, the DSO Ringerikskraft provided only a volumetric cost at $0.02425 \frac{EUR}{kWh}$. With the proposed MP capacity-based grid tariff from NVE, accumulated for a monthly period \cite{Norgesvassdrags-ogenergidirektorat2017ForslagNettvirksomhet}, the volumetric cost would be at $0.00625 \frac{EUR}{kWh}$, and a monthly capacity-based cost at $7.2075 \frac{EUR}{kW_{peak}}$.


\subsection{Scenario Generation}

The HEMS together with the SDP algorithm allows the possibility for multiple input data to be uncertain in the period of operation. To limit the range of uncertainty, the work here considers uncertainty within the EV behavior, outdoor temperature and solar irradiation. Information such as electricity price and electric-specific demand is considered deterministic. 

In total, 9 scenarios per day have been generated for this case study, where 3 scenarios have been generated from both weather effects and EV behavior. It is assumed that the sources are independent of each other, resulting in 9 combinations. The stochastic nature of the input data is based on a normal distribution with the mean and standard deviation as the discrete scenarios, giving a probability distribution at $\rho_{\mu} = 68.2 \%, \rho_{\sigma} = 15.9 \% $ for each source. 

The normal distribution of EV behavior for a weekday is based on \cite{LakshmananValueScenario}, with an expected departure/arrival time from 9 AM to 5 PM, and a standard deviation of 90 minutes (as this work considers an hourly time step, the standard deviation is rounded up to 2 hours). 
Moreover, authors in \cite{Sadeghianpourhamami2018QuantitiveApproach} show that the expected arrival time when charging near home does not change dramatically between weekdays and weekends. Thus, we assume the same departure/arrival time pattern for the weekend.

Data for both the outdoor temperature and solar irradiation have been obtained from Rygge weather station in South-east Norway \cite{NedlastingLMT}. Hourly data from 2014-2019 for the month of January have been used to create hourly normal distributions, to generate three discrete scenarios. 

\subsection{Model Cases}

 The scope of this work will be to investigate the flexibility contribution each flexible asset can provide, both for MPGT and energy-based grid tariff (EBGT) structures, both given $MPGT$ and $EBGT$ as case names, respectively. This extends the analysis in \cite{EmilThorvaldsen2020RepresentingProgramming} by investigating the long-term value each flexible asset offers individually instead of combined. For each case with a specified flexible asset, the other assets will have a passive behavior. The passive manner for each of the assets is the following: The BESS is turned off, the EV charging will charge to max capacity at 90\% whenever it can and has an initial condition of $E_0^{EV} = 22.6$ kWh during transition. The space heating will maintain a constant indoor temperature at $T_t^{in} = 22$ $^\circ C \: \forall t$. Moreover, the flexible assets will be investigated for several input parameter values to analyze the change of impact, as showcased in Table \ref{tab:parameter_input}.
 
 The analysis will be carried out by first generating the EFCCs for each case using the SDP algorithm in Alg. \ref{Alg_1}. The state variable subset $\mathcal{S}_{P,g}$ will consist of 100 discrete initial values of $p^{imp}$, ranging from 0-10 $\frac{kWh}{h}$. This leads to 27 900 unique decomposed problems to solve per case. After this, the value of the EFCCs will be analyzed through a simulation phase, where each day is run sequentially for the whole month to see the overall economical performance, where the peak import level will be carried over between stages and the final grid tariff cost is set at the end. To account for the range of uncertainty in the input data, the simulation phase analyze this month 1000 times with different stochastic scenario realizations based on their scenario probability. The sequential coupling between two stages is done by randomly drawing a future scenario for the next stage transition, where the odds for each scenario is based on their probability. This leads to multiple scenario combinations per month.

\begin{table}[]
\centering
\caption{Flexible asset parameters for the different cases. Values in bold are default values when considering the asset as passive.}
\begin{tabular}{ccc}
\hline
\textbf{Component}        & \textbf{Parameter(s)}            & \textbf{Cases}                                                     \\ \hline
\textbf{Battery energy storage system}          & \textbf{$E^{B,Cap}$}             & \textbf{5 kWh}, 10 kWh                            \\ \hline
\textbf{Space heating}    & \textbf{$T^{in,min},T^{in,max}$} & \textbf{{[}20, 24{]}}, {[}21, 23{]} \\ \hline
\textbf{Electric Vehicle} & \textbf{$\dot{E}^{Max}$}         & 2.3 kW, \textbf{3.7 kW}                    \\ \hline
\end{tabular}\label{tab:parameter_input}
\end{table}

%% file: chapters/Results_and_discussions.tex
\section{Results \& Discussion} \label{discussion}

The SDP algorithm showcased in Algorithm \ref{Alg_1} creates EFCC from the last day of the month, and by working backwards to create an accurate EFCC for the first stage, that presents the future cost associated with the peak import state variable. 
To obtain an overview of the capabilities and potential given by the generation of the EFCCs, case $MPGT$ will be analyzed first, with the economic performance and EFCCs presented in Section \ref{eco_performance}. Furthermore, $MPGT$ will be compared to $EBGT$ in Section \ref{energy_power_discussion}.

\subsection{Long-term price signal performance $MPGT$} \label{eco_performance}

\subsubsection{Economic performance}

\begin{figure}[h!]
\centering
\includegraphics[width=0.8\textwidth]{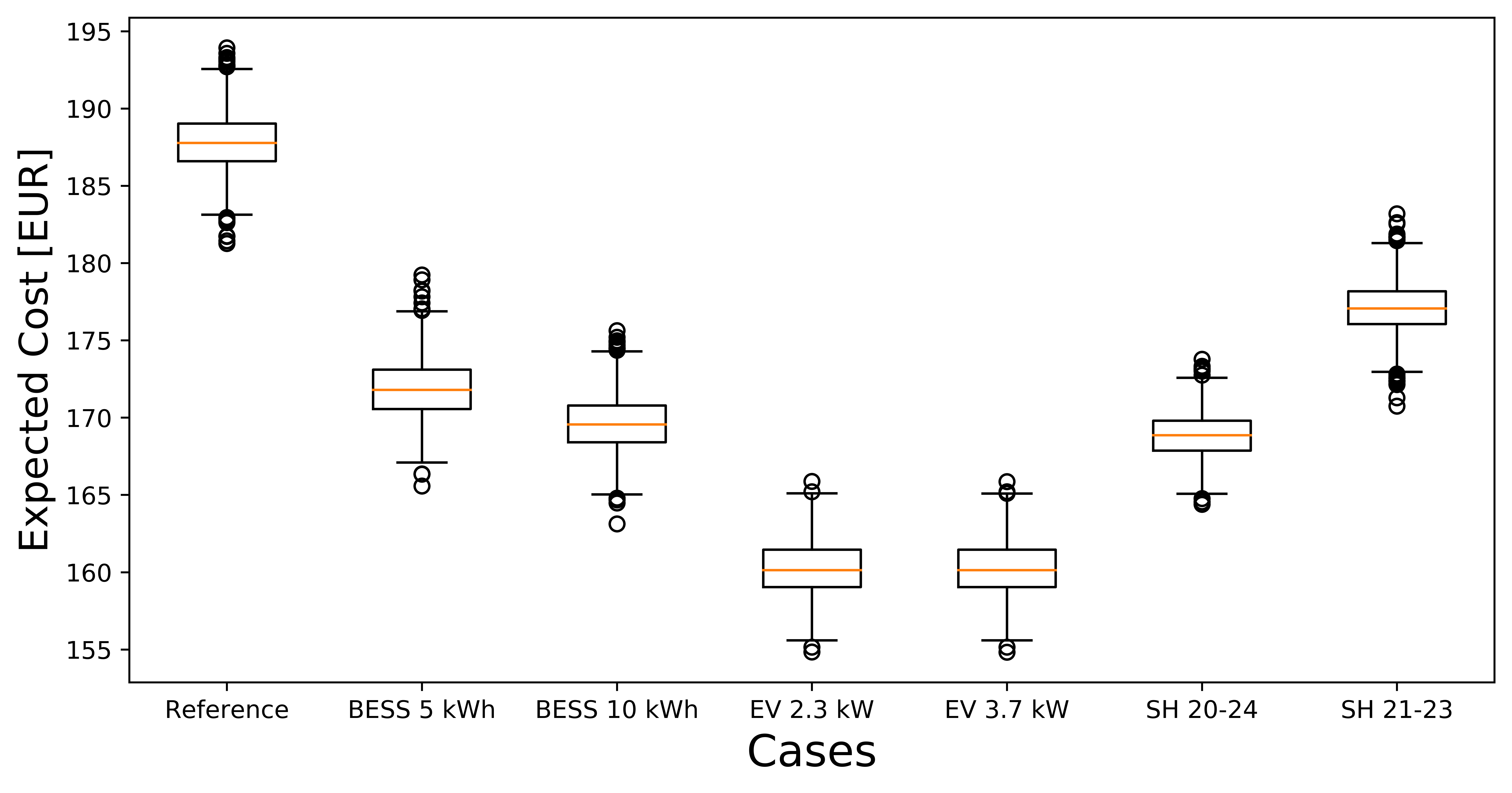}
\caption{Boxplot of expected monthly total cost for the different cases for scheme $MPGT$.}
\label{fig:Cost_RTP}
\centering
\end{figure}

As the case study presented here involves uncertainty, the economic performance will vary based on the sequence of scenario realizations. To capture the tendencies and dispersion of the data, the expected total cost for the month based on the different cases is plotted as a boxplot in Fig. \ref{fig:Cost_RTP}. Within the boxplot, the median value is indicated by the orange line, whereas the box specifies the interquartile range (IQR) of the results. The lines outside the boxes are whiskers, representing the 1.5*IQR, while the outliers show the few results that are outside of the 1.5*IQR.

Fig. \ref{fig:Cost_RTP} illustrates how every flexible asset manages to decrease the cost in comparison to the reference case. Both EV charging cases have the highest impact with an expected cost decrease of 14.6\% compared to the reference case, with $SH_{20,24}$ being the second best with a 10.1\% decrease, in which the latter had a noticeable positive cost change compared to the other SH case. The 10 kWh BESS had a 9.7\% cost reduction, coming up close to the best SH case. The range of possibilities in total cost due to the uncertainty, shows that the cost has a specific span of possible results, in which all cases have a similar IQR and whisker range.

\begin{figure}[h!]
\centering
\includegraphics[width=0.8\textwidth]{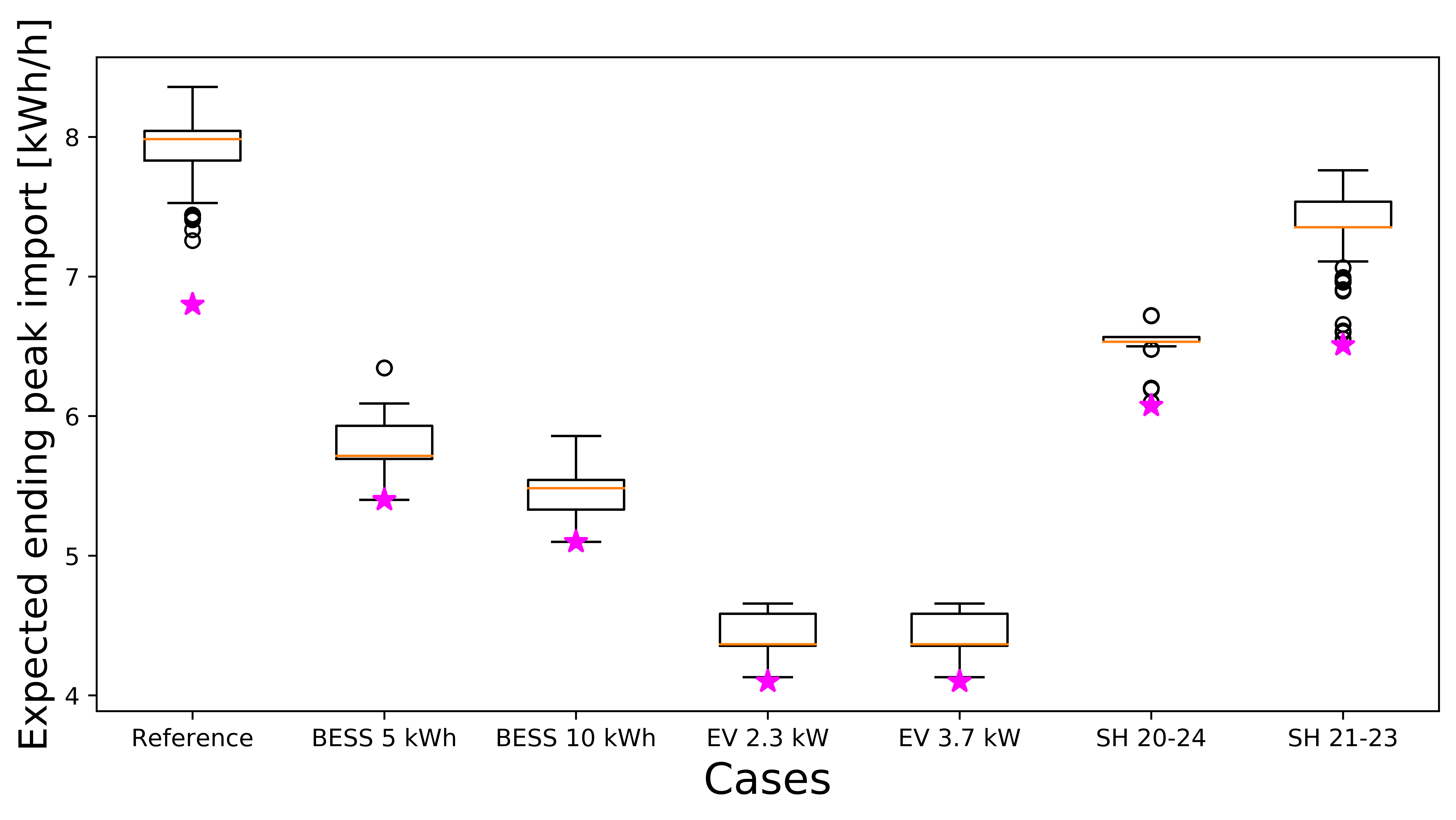}
\caption{Overview of expected ending peak for the different cases for $MPGT$. The star is the initial peak import level from day 1.}
\label{fig:peak_RTP}
\centering
\end{figure}

As the economic performance in Fig. \ref{fig:Cost_RTP} is affected by the MPGT in addition to RPT, the ending peak level for each case is showcased in Fig. \ref{fig:peak_RTP} to illustrate the MPGT impact. All cases achieve a reduction in ending peak level, although the amount depends on which flexible asset is activated.

For the BESS cases, they can charge/discharge strategically to reduce the peak import as shown in Fig. \ref{fig:peak_RTP}. The battery is not a natural load within the building, making the value of flexibility unique compared to the rest with no risk of any rebound effect. With an inverter of 2.5 kW, the expected peak reduction from the reference case is at 2.20 kW and 2.50 kW for the 5 and 10 KWh BESS, respectively. Thus, there is a correlation to not only power capacity, but also storage capacity, indicating that the peaks last for multiple hours. However, the 5 kWh BESS manages to reduce a relatively large proportion of the peak compared to the 10 kWh BESS. 
For the EV charger cases, the flexible contribution provides the highest peak import reduction, with an expected peak import cut by 3.54 kW. As the peak result is similar for both capacities, with the reference case utilizing a 3.7 kW charger, this showcases how significant an impact the passive EV charging strategy has, and how vital flexible EV charging is. As it is possible to charge for 15 hours per day, there is a high range of flexibility to choose from, making the charging capacity minor if flexible. However, as the EV is uni-directional, it can only decrease peak from its own demand. 

For space heating, the performance improves as the indoor temperature boundary increases. The heater is used diligently during this cold winter month to keep the temperature in check. 
Case $SH_{21,23}$ has the lowest total cost reduction and peak import reduction compared to the other cases. The cost decrease is derived from both the peak reduction, and by the utilization of flexibility for RTP adjustments, whereas the latter seems to have the most impact. However, by looking at the peak import, the capability of load shifting in peak hours for $SH_{21,23}$ is deemed less critical in the long run. If the boundary was increased further to 20-24 $^\circ C$, the total cost will be reduced as well as the peak is declined further. The deviation between the two SH cases shows how the temperature boundary affects the load shifting capability, where a higher boundary gives more capacity of pre-heating the interior before peak hours while avoiding high RTP hours as a rebound effect, allowing a longer idle period. Based on the differences, the 21-23 $^\circ C$ boundary, with the limitations in continuous idle hours, shows it to be more cost-optimal to increase peak to decrease RTP costs. Both cases find the threshold for peak power that gives the optimal balance of MPGT and RTP costs seen together, and continuously considers the marginal MPGT increase vs decreasing RTP costs for each scheduling day.

An interesting takeaway from the SH results, is the peak import for both cases, which from Fig. \ref{fig:peak_RTP}
is likely to change during scheduling from the initial point. Case $SH_{21,23}$ has a much wider range of ending peak import than the other. 
In $SH_{20,24}$, the first and third quartiles, as well as the whiskers, are very tight for the peak import, indicating that the added flexibility gives it more room to reach the same ending peak import, almost regardless of the scenarios realized. Moreover, the deviation from the initial point is smaller than for $SH_{21,23}$, which shows the latter case is more prone to scenario realizations, as that would affect RTP costs. For $SH_{20,24}$, these results show how adjustable SH can be regarding peak import, since as the total cost has a wider spread, it keeps the peak import stable, and instead utilizes load shifting which increases RTP costs, but which is deemed as less costly than the increased grid tariff cost. It is worth mentioning that when increasing the temperature boundary beyond $SH_{20,24}$, the changes are marginal, indicating that this boundary gives the most valuable flexibility capability increase for this case, fitting well with the observation from \cite{Berge2016PerceivedBuildings}. 

\begin{figure}[h!]
\centering
\includegraphics[width=0.8\textwidth]{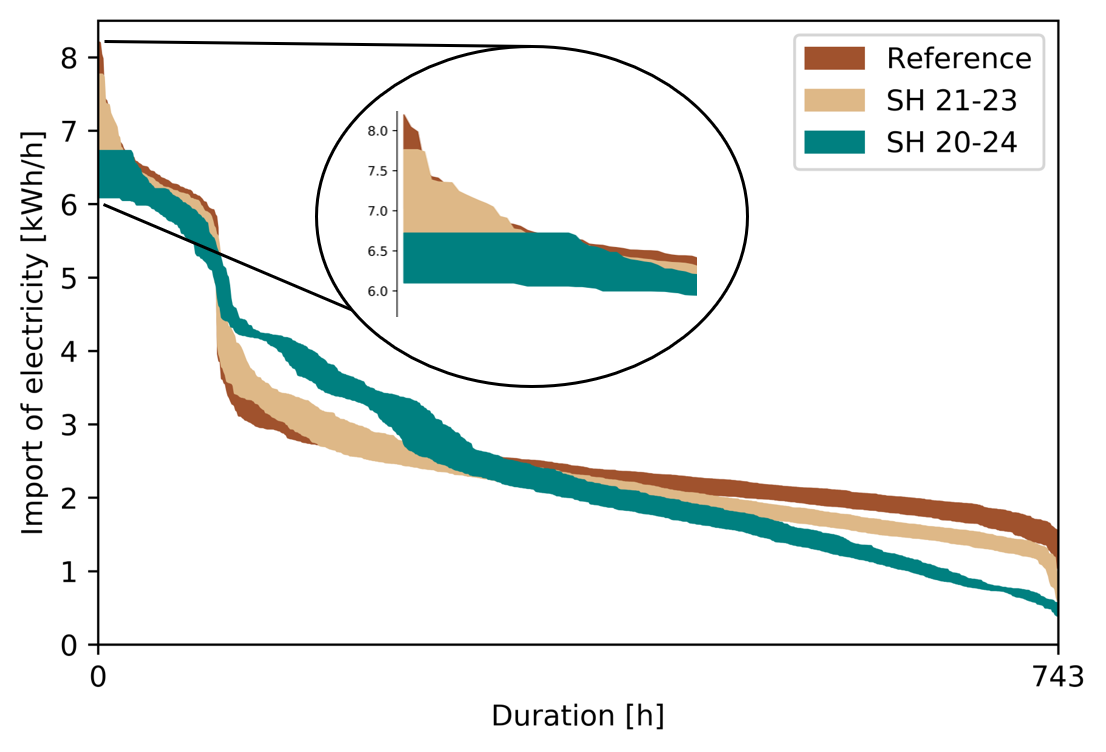}
\caption{Duration curve for the specified space heating cases over a month, with the uncertainty of scenarios included.}
\label{fig:SH_duration}
\centering
\end{figure}

Moreover, Fig. \ref{fig:SH_duration} illustrates the import duration curve for the reference and the two SH cases. The most prominent behavior is how case $SH_{20,24}$ manages to cut the peak through the use of flexibility, whereas case $SH_{21,23}$ has higher peak demand as some critical scenarios makes a higher peak more beneficial than increasing RTP cost as a rebound effect, overlapping with the reference case at the peaks.

\subsubsection{Marginal expected future cost curves} \label{EFCC_discussion}

\begin{figure}[h!]
\centering
\includegraphics[width=0.8\textwidth]{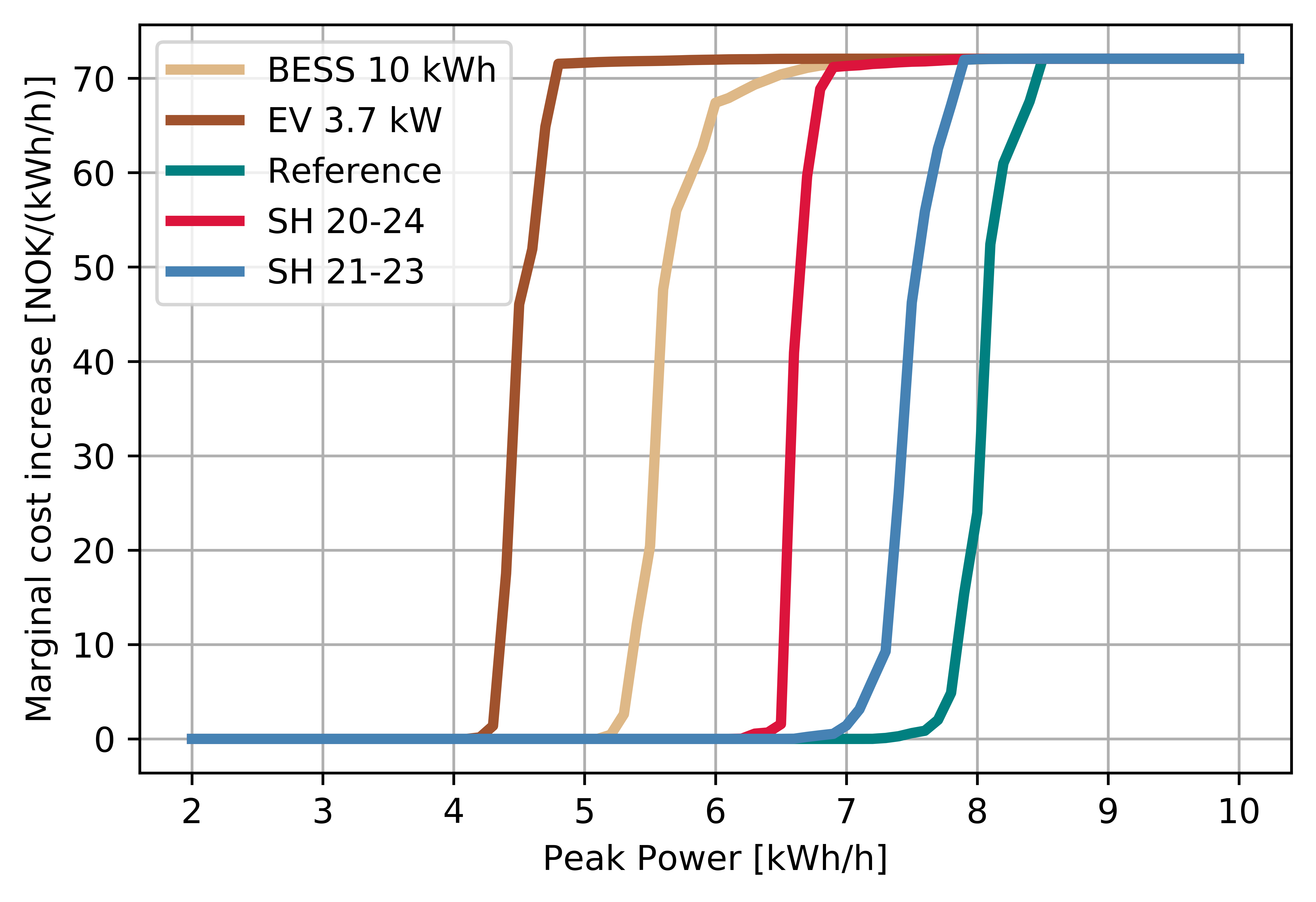}
\caption{Plot of the MEFCC for the different cases on day 1.}
\label{fig:EFCC_all}
\centering
\end{figure}

To demonstrate how the EFCC changes, the different cases are shown in Fig. \ref{fig:EFCC_all}, plotted as marginal EFCC (MEFCC) to make a better comparison of the marginal change based on the peak import. 

 In general, the curves can be divided into three main areas. For marginal cost increase at 0, the initial peak power will not be economically suitable to stay at this level, either because of the lack of flexibility potential, or through costly use of flexibility that is not worth it. Likewise, for the part where the marginal cost increase is equal to the grid tariff cost at 7.25 $\frac{EUR}{kW_{peak}}$, there is no future cost-benefit other than the increasing grid tariff at the end. Between these two areas, the curve is increasing in cost from 0 to the grid tariff cost, which shows the value of enabling more flexibility from peak import levels. This area is influenced by the reaction of flexible assets towards peak levels, and the uncertainty of the problem. In comparison to the flexible assets, all assets manage to shift this intersection towards the left, thus highlighting the capability to cut peak import.

The 10 kWh BESS decreases the peak import by charging during non-peak hours, and discharge during critical hours, cutting the peak by around 2.5 kW compared to the reference case. 
The curve has a marginal cost increase with increasing peak power, which shifts to the right to a higher degree than all the other cases, due to scenario realizations and efficiency loss savings with less usage. By increasing the inverter capacity, the peak shaving could be more substantial, however this could lead to the situation with the 5 kWh BESS where the power-energy rating would come into question for the behavior of the MEFCC. The MEFCC of EV charger has the lowest peak import level change compared to the reference case. The EV charger cuts the peak level substantially by being flexible in charging, which is primarily due to the significant impact the charger gives if being passive, such as in the reference case. Finally, the MEFCC for both $SH_{20,24}$ and $SH_{21,23}$ have different behavior. As stated in Section \ref{eco_performance}, the temperature boundary affected the peak demand level and the variance significantly, which is confirmed by the shape of the curves. The $SH_{21,23}$ curve has the smallest peak level decrease from the reference case, and is scenario sensitive with the shift towards the right, whereas the $SH_{20,24}$ curve shows little scenario impact due to the almost vertical cost increase.

\subsection{The economic comparison of $MPGT$ and $EBGT$} \label{energy_power_discussion}

\begin{table*}[t]
\centering
\caption{The economic performance of all sensitivity cases on both electricity schemes. Values given in EUR}
\begin{tabular}{ccc}
\hline
Case \textbackslash{} Scheme & EBGT & MPGT \\ \hline
Reference                                     & 167.8                & 187.79   \\ \hline
BESS 5 kWh                                  & 167.38                & 171.88                         \\ \hline
BESS 10 kWh                               & 167.17                & 169.61                             \\ \hline
EV 2.3 kW                                 & 165.8                & 160.29                           \\ \hline
EV 3.7 kW                                   & 165.8                & 160.28    \\ \hline
$SH_{20,24}$                                   & 156.65               & 168.88          \\ \hline
$SH_{21,23}$                                   & 159.41                & 177.14         \\ \hline

\end{tabular} \label{tab:Eco_perf}
\end{table*}

The economic performances of all cases are found in Table \ref{tab:Eco_perf}, showing the average total cost.

By analyzing each scheme separately, we can see how the expected total cost decreases compared to the reference case when the flexible assets are enabled to be controlled. For $EBGT$, the cost decrease is tied to the load-shifting capability to utilize the RTP deviations. For this electricity scheme, it is mostly defined for the SH cases, as the high use for heating during winter months leads to increased consumption, which gives higher cost saving potential if controllable. 

When comparing the schemes, it is evident that the average cost of $MPGT$ is higher than for $EBGT$, except when activating smart EV control, which achieves the lowest expected total cost with the $MPGT$ scheme. In addition, $SH_{20,24}$ manages to achieve the next best result in $MPGT$ despite having a higher ending peak than the batteries, due to load shifting to balance peak power and RTP costs. The utilization of load shifting for RTP is also seen in the $EBGT$ scheme, demonstrating the high potential in cost reduction by utilizing flexibility for both peak power and RTP. The 10 kWh BESS has a marginal contribution to RTP alone in the $EBGT$ scheme, but provides a much higher cost reduction within the $MPGT$ scheme by reducing the peak import.


However, what these two electricity schemes provide is an overview of the future value of the flexibility that the assets can offer. With the $EBGT$ scheme that only benefits RTP deviations, the cost decrease of smart control is generally lower than with the MPGT, compared to the reference cases. Given that the grid tariff is needed to reflect the cost of operating more accurately, more time-based cost deviations might occur or be presented, which as presented here can provide much more cost savings than if we remain with the default behavior. Flexible assets are valuable even during current operation, however their potential can only increase in the future, both in regards to short- and long-term price signals as we have presented here.

%% file: chapters/Conclusion.tex
\section{Conclusion} \label{Conclusion}

We have presented a model that aims to illustrate the expected future cost for a building operation model with a long-term price signal. 
The future cost is based on the expected succeeding cost for operation for a building and the long-term price signal, in this work being a measured-peak grid tariff. This model was applied to a Norwegian household. The primary goal was to analyze the contribution from a battery energy storage system, a smart Electric vehicle charger, and controllable space heating individually with varying input parameters, to see how they react to cope with the long-term price signal. The performance was compared to two electricity schemes with and without the long-term price signal. 

The results from the generation of expected future cost curves showed that all flexible assets contribute to lower the peak import level compared to a reference case, but revealed that their flexibility characteristics affect the long-term performance. The generation of expected future cost curves enables us to represent the future impact of current short-term decision-making, which can provide more cost-optimal flexibility usage over the total period. The results here have shown the contribution of each flexible asset individually, whereby this also suggests that further investigation of their contribution together and for other price signals could enable a broader understanding of their capabilities. 

\section{Acknowledgment} \label{Acknowledgment}

This work was funded and supported by the Research Council of Norway (Grant Number: 257626/257660/E20) and several partners through FME ZEN and FME CINELDI. The authors gratefully acknowledge the financial support from the Research Council of Norway and all partners in CINELDI and ZEN.